\documentclass[12pt]{article}
\usepackage{amsmath,amssymb,pb-diagram}

\def \({\left(}
\def \){\right)}
\def \1{{\bf{1}}}

\def \al{\alpha}

\def \bs{\backslash}

\def \C{{\mathbb C}}

\def \CF{{\cal F}}

\def \CM{{\cal M}}

\def \CP{{\cal P}}

\def \CO{{\cal O}}
\def \CS{{\cal S}}

\def \diag{{\rm diag}}

\def \ds{\displaystyle}

\def \eps{\varepsilon}
\def \eqn{\begin{eqnarray*}}
\def \endeqn{\end{eqnarray*}}

\def \ga{\gamma}
\def \Ga{\Gamma}

\def \H{{\mathbb H}}

\def \Im{{\rm Im}}

\def \La{\Lambda}

\def \N{\mathbb N}

\def \ph{\varphi}
\def \prf{\vspace{5pt}\noindent{\bf Proof: }}
\def \P{{\mathbb P}}

\def \PSL{{\rm PSL}}

\def \qed{\ifmmode\eqno $\square$\else\noproof\vskip 
			12pt plus 3pt minus 9pt \fi}
 \def\noproof{{\unskip\nobreak\hfill\penalty50\hskip2em\hbox{}%
     \nobreak\hfill $\square$\parfillskip=0pt%
     \finalhyphendemerits=0\par}}

\def \R{{\mathbb R}}
\def \Re{{\rm Re \hspace{1pt}}}

\def \setminus{
	\begin{picture}(18,10)\put(4,6){\line(2,-1){10}}\end{picture}}
\def \sign{\mathrm{sign}}

\def \SL{{\rm SL}}

\def \tr{{\hspace{1pt}\mathrm{tr}\hspace{2pt}}}

\def \z{\zeta}
\def \Z{\mathbb Z}
\def \={\ =\ }

\newtheorem{theorem}{Theorem}[section]

\newtheorem{lemma}[theorem]{Lemma}
\newtheorem{remark}[theorem]{Remark}

\newtheorem{proposition}[theorem]{Proposition}
\newcommand{\norm}[1]{\left|\hspace{-1pt}\left| #1\right|\hspace{-1pt}\right|}
\renewcommand{\matrix}[4]{\begin{pmatrix}#1 & #2 \\ #3 & #4 \end{pmatrix}}
\renewcommand{\sp}[1]{\left\langle #1 \right\rangle}

\begin{document}
\pagestyle{myheadings} \markright{THE LEWIS CORRESPONDENCE...}

\title{The Lewis Correspondence for submodular groups}
\author{Anton Deitmar \& Joachim Hilgert}
\date{March 25, 2004}
\maketitle
\tableofcontents

\newpage
\section{Introduction}

Let $G$ denote the group $\PSL_2(\R)=\SL_2(\R)/\{\pm 1\}$. Let
$\Ga=\PSL_2(\Z)$ be the modular group.
A \emph{submodular group} is a subgroup of  $\Ga$ of finite
index. It is the aim of this note to extend the Lewis Correspondence
\cite{Lewis, Lewis-Zagier1, Lewis-Zagier2} from $\Ga$ to submodular groups.
Since any submodular group $\La$ contains a submodular subgroup which is
normal in $\Ga$ we will first assume that $\La$ is normal and only later
move from $\Lambda$ an arbitrary subgroup containing $\Lambda$. Let
$\H^+=\{ x+iy\in\C: y>0\}$ be the upper half plane in $\C$. The group $G$ 
acts on $\H^+$ by linear fractions,
$\ds\matrix abcd \cdot z=\frac{az+b}{cz+d}$.
This action preserves the hyperbolic geometry given by the
Riemannian metric $\frac1 {y^2}(dx^2+dy^2)$  so it commutes with the 
hyperbolic Laplace operator
$\Delta=-y^2\left((\frac\partial{\partial x})^2
+(\frac\partial{\partial y})^2\right)$
and preserves the hyperbolic volume form $dx\, dy/y^2$.
A \emph{Maa\ss\  form} for $\La$ is a function
$f\in L^2(\La\bs\H^+)$ which is an eigenfunction of $\Delta$.
We also define $\H^-$ to be the lower half plane in $\C$.
The Lewis Correspondence
attaches a certain ``period function'' to a given Maa\ss\ form
for $\Ga$. To extend it to
$\La\ne \Ga$ we have to start with Maa\ss\ forms for
$\La$. These form a module under the
finite group $\Ga/\La$ under the regular representation
and so Maa\ss\
forms for $\La$ are related to Maa\ss\ forms for
$\Ga$ twisted by a finite dimensional representation $(\eta, V_\eta)$ by
the following mechanism:

Let $W$ be a $\C[\Ga]$-module, which is finite
dimensional as $\C$-vector space and trivially acted upon by $\La$. Under the
action of the finite group $\Ga/\La$ the module $W$ decomposes into isotypic
components,
\begin{equation}\label{isotypegeneral}
W\=\bigoplus_{\eta\in\widehat{\Ga/\La}}W(\eta),
\end{equation}
where $\widehat{\Ga/\La}$ denotes
the set of isomorphism classes of
irreducible unitary representations of $\Ga/\La$, i.e., the
unitary dual of
this finite group. For $\eta\in\widehat{\Ga/\La}$ let $\breve\eta$ denote its
dual representation. There is a natural isomorphism
\begin{equation}\label{isotypisom}
\mathrm{ev}\colon(W\otimes\eta)^\Ga\otimes\breve\eta\ \to\ W(\breve\eta)
\end{equation}
given by
$\mathrm{ev}(\sum_j (w_j\otimes\al_j)\otimes\beta):=\sum_j\sp{\al_j,\beta}w_j$.
On the other hand, the inclusion $W(\breve\eta)\subset W(\eta)$ induces an isomorphism $(W\otimes\eta)^\Ga\cong(W(\breve\eta)\otimes\eta)^\Ga$ and the projection map $\Pr$ from $W\otimes\eta$ to $(W\otimes\eta)^\Ga$ is explicitly given by
$$
\Pr(w\otimes\al)\=\frac 1{|\Ga/\La|}\sum_{\ga : \Ga/\La}\ga .w\otimes\ga.\al.
$$
Finally, elementary character theory shows that the canonical projection
$\CP_{\breve\eta}\colon W\to W(\breve\eta)$
given by the decomposition (\ref{isotypegeneral}) equals
\begin{equation}\label{isotypproj}
\CP_{\breve\eta}w
=\frac{d_\eta}{|\Gamma/\Lambda|}\sum_{\gamma\in \Gamma/\Lambda}
  \tr\eta(\gamma)\ (\gamma\cdot w),
 \end{equation}
where $d_\eta$ ist the degree of $\eta$ and $\breve\eta$.
Here we have used the convention that we
write the space of a representation with the same symbol as the
representation itself. Occasionally, to put emphasis on the space rather
than the representation, we will also write $V_\eta$ for the representation
space of $\eta$. In order to describe
$W$ we decompose it into isotypic components and each such component is
described by $(W\otimes\eta)^\Ga$. 
We will in particular apply this to the space of Maa\ss\ forms for $\La$ with a given Laplace eigenvalue.
But we also can retrieve Maa\ss\ forms of an arbitrary submodular group $\Sigma$.
For this let $\La\subset\Sigma\subset\Ga$ be a submodular group which is normal in $\Ga$, and let $W$ be the space of $\La$-Maa\ss\ forms.
Then
$$
W\ \cong\ \bigoplus_{\eta\in\widehat{\La/\Ga}}(W\otimes\eta)^\Ga\otimes\breve\eta.
$$
The space of $\Sigma$-Maa\ss\ forms is just the space of $\Sigma$-invariants herein, i.e., the space
$$
W^\Sigma\ \cong\ \bigoplus_{\eta\in\widehat{\La/\Ga}}(W\otimes\eta)^\Ga\otimes\breve\eta^\Sigma.
$$
So $W^\Sigma$ is described by the spaces $(W\otimes\eta)^\Ga$ and the dimensions $\dim(\breve\eta^\Sigma$ for $\eta\in\widehat{\La/\Ga}$.
This applies in particular to the congruence subgroups $\La=\Ga(N)$ and $\Sigma=\Ga_0(N)$.
 So we fix an irreducible representation $\eta$ of
$\Ga$ with finite image.

We fix the following notation for the canonical 
generators of $\Ga$:
$$
S\=\pm\matrix 01{-1}0,\quad {\rm and} \quad T\=\pm\matrix 1101.
$$
Then $S^2=\1=(ST)^3$, and $T$ is of infinite order.
Let $\CF_\eta$ be the space of
holomorphic functions $f\colon
\C\setminus\R\to V_\eta$ with
\begin{eqnarray}
\label{Fetaa} f(z+1)&=&\eta(T) f(z),\\
\label{Fetab} f(z)&=&O(1)\quad  \text{ as }|\Im(z)|\to\infty,\\
\label{Fetac} 0&=&f(i\infty)+f(-i\infty).
\end{eqnarray}
The last condition needs explaining. Since $\eta$ has finite image, there is
a smallest $N:=N_\eta\in\N$ such that $\eta(T^N)$ equals the identity.
It follows that $f$ has a Fourier expansion
$$
f(z)\=\sum_{k\in\Z} e^{2\pi i\frac kN z} v_k^+,\qquad v_k^+\in V_\eta,
$$
in $\H^+$ and similarly with $v_k^-\in V_\eta$ in $\H^-$. Condition
(\ref{Fetab}) leads to $v_{-k}^+=v_k^-=0$ for every $k\in\N$. Thus the limits do exist
and satisfy $f(i\infty)=v_0^+$ and $f(-i\infty)=v_0^-$.

Consider the space $\CF_{\nu,\eta}$ of all $f\in\CF_\eta$
for which the map
\begin{equation}\label{Fnueta}
z\mapsto f(z)-z^{-2\nu- 1}\,\eta(S) f\left(\frac{-1}z\right)
\end{equation}
extends holomorphically to
$\C\setminus (-\infty,0]$ and 
the space $\Psi_{\nu,\eta}$ of all holomorphic functions $\psi\colon
\C\setminus(-\infty,0]\to V_\eta$ satisfying
\begin{equation}\label{Psinuetaa}
\eta(T)\psi(z)=
\psi(z+1)+(z+1)^{-2\nu-1}\eta(ST^{-1})\psi
\left(\frac z{z+1}\right)
\end{equation}
and
\begin{equation}\label{Psinuetab}
\begin{array}{l}
0\ =\   \ds e^{-\pi i \nu}\lim_{\Im(z)\to\infty}
             \left(\psi(z)+z^{-2\nu-1}\eta(S)\psi\left(\frac {-1}z\right)\right)+\\
\phantom{0\ =\ } \ds  +e^{\pi i\nu}\lim_{\Im(z)\to-\infty}
             \left(\psi(z)+z^{-2\nu-1}\eta(S)\psi\left(\frac{-1}z\right)\right),\\
   \end{array}
\end{equation}
where both limits exist. We call   (\ref{Psinuetaa}) the
\emph {Lewis equation}.

Let $\pi_\nu$ be the principal series representation of $G$ associated with the 
parameter $\nu\in \C$ and $\pi_\nu^{-\omega}$ the corresponding
space of hyperfunction vectors.
As a crucial tool we will use the space 
\begin{equation}\label{etaautomhypfcts}
A_{\nu,\eta}^{-\omega}=(\pi_\nu^{-\omega}\otimes\eta)^\Ga=
H^0(\Ga,\pi_\nu^{-\omega}\otimes\eta)
\end{equation}
and call it the space of \emph{$\eta$-automorphic hyperfunctions}.

Generalizing results of Bruggeman (see \cite{Brug}, Prop.~2.1 and Prop.~2.3), 
we will show in Proposition \ref{1.1}
that there is a linear isomorphism
$A_{\nu,\eta}^{-\omega} \to \CF_{\nu,\eta}$
and (using this) establish in Proposition \ref{1.2} a linear map
$$
B\colon A_{\nu,\eta}^{-\omega}\to \Psi_{\nu,\eta}, 
$$
which we call
the \emph{Bruggeman transform}. It turns out to be bijective unless
$\nu\in \frac 12+\Z$.

Recall that a \emph{Maa\ss\ wave form} for a subgroup $\La$ of $\Gamma$ (not
necessarily
normal) and parameter $\nu\in\C$ is a
function $u$ on $\H^+$ which is twice continuously differentiable and
satisfies
\begin{eqnarray}
\label{Mnua} u(\ga z)&=&u(z) \quad \text{ for every }\ga\in\Ga,\\
\label{Mnub} \infty &>& \ds\int_{\Ga\bs\H^+}|u(z)|^2\, dz,\\
\label{Mnuc} \Delta u&=&\({\textstyle\frac 14 -\nu^2}\) u.
\end{eqnarray}
By the regularity of solutions of elliptic differential equations the
last condition implies that $u$ is real analytic.
Let $\CM_\nu=\CM_\nu^\La$ be the space of all Maa\ss\ wave forms for $\Lambda$.

If $\Lambda$ is normal of finite index in $\Gamma$ the
finite
group $\Ga/\La$ acts on this space, and as in (\ref{isotypegeneral})
we get an isotypic decomposition,
\begin{equation}\label{MaassLambdaisotypes}
\CM_\nu\=\bigoplus_{\eta\in\widehat{\Ga/\La}}\CM_\nu(\eta).
\end{equation}
and for each $\eta$,
$$
\CM_\nu(\breve\eta)\ \cong\ V_\eta^*\otimes\( V_\eta\otimes\CM_\nu\)^{\Ga/\La}.
$$
We  set
$\CM_{\nu,\eta}$ equal to $\( V_\eta\otimes\CM_\nu\)^{\Ga/\La}$. Then $\CM_{\nu,\eta}$
can be viewed as the space
of all functions $u\colon\H^+\to V_\eta$ which are twice continuously
differentiable and satisfy
\begin{eqnarray}
\label{Mnuetaa} u(\ga z)&=&\eta(\ga)u(z) \quad \text{ for every }\ga\in\Ga,\\
\label{Mnuetab} \infty&>&\ds\int_{\Ga\bs\H^+}\norm{u(z)}^2\, dz, \\
\label{Mnuetac} \Delta u&=&\({\textstyle\frac 14 -\nu^2}\) u.
\end{eqnarray}

We define the space $\CS_{\nu,\eta}$ of \emph{Maa\ss\ cusp forms} to be the
space of all $u\in\CM_{\nu,\eta}$ such that
\begin{equation}\label{Snueta}
\int_0^N u(z+t)\, dt\= 0
\end{equation}
for every $z\in\H^+$. Here, as before, $N$ is the order of $\eta(T)$, so that in particular
$\eta(T)^N=\1$ and  $u(z+N)=u(z)$.

For $\Re \nu>-\frac 12$ consider the space $\Psi_{\nu,\eta}^o$
of all $\psi\in \Psi_{\nu,\eta}$  satisfying
\begin{equation}\label{Lewisasymcomplex}
\psi(z)\ =\ O(\min\{1,|z|^{-C}\}) \quad \text{ for }  z\in \C\setminus(-\infty,0],
\end{equation}
for some $0<C<2\Re \nu+1$.
We call  the elements of $\Psi_{\nu,\eta}^o$ \emph{period functions}.
In Lemma \ref{Lewistrafolemma} we establish for 
$\Re \nu>-\frac 12$ two linear maps
$\CS_{\nu,\eta}\to \CF_{\nu,\eta}$ and
$\CS_{\nu,\eta}\to \Psi_{\nu,\eta}^o$.

Let $\tilde\CM_{\nu,\eta}$ be the
space of all functions $u$ satisfying only (\ref{Mnuetaa}) and (\ref{Mnuetac}).
So there is no growth restriction on elements of $\tilde\CM_{\nu,\eta}$.
For an automorphic hyperfunction $\al\in A_{\nu,\eta}^{-\omega}$ we consider the
function $u\colon G\to V_\eta$ given by
$$
u(g) \ :=\  \sp{\pi_{-\nu}(g)\ph_0,\al}.
$$
Then $u$ is right $K$-invariant, hence can be viewed as a function on $\H^+$.
As such it lies in $\tilde\CM_{\nu,\eta}$ since $\al$ is $\Ga$-equivariant and
the Casimir operator on $G$, which induces $\Delta$, is scalar on $\pi_\nu$
with eigenvalue $\frac 14-\nu^2$.
The transform $P\colon \al\mapsto u$ is called the \emph{Poisson transform}. It
follows from \cite{Schlichtkrull}, Theorem 5.4.3, that the Poisson transform
\begin{equation}\label{Poissonisom}
P\colon A_{\nu,\eta}^{-\omega}\ \to\ \tilde\CM_{\nu,\eta}
\end{equation}
is an isomorphism for $\nu\not\in \frac 12+\Z$.

For $\nu\not\in \frac 12+\Z$ we finally define the \emph{Lewis transform} as the
map $L\colon \CM_{\nu,\eta}\to \Psi_{\nu,\eta}$,
given by
\begin{equation}\label{Lewistrafo}
L\ :=\  B\circ P^{-1}.
\end{equation}
Our first main result (see Theorem \ref{Lewistrafobijektiv})
is a generalization of \cite{Lewis-Zagier2}, Thm.~1.1, and says
that the Lewis transform for $\nu\not\in \frac 12+\Z$ and 
$\Re \nu> -\frac 12$
is a linear isomorphism between the space of Maa\ss\ cusp
forms $\CS_{\nu,\eta}$ and the space $\Psi_{\nu,\eta}^o$ 
of period functions.

A holomorphic function on $\C\setminus (-\infty,0]$ is uniquely determined by its
values in $\R^+:=(0,\infty)$. Thus, in principle, it is possible to describe the period
functions as a space of real analytic functions on the positive halfline. Following
ideas from \cite{Lewis-Zagier2}, Chap. III, in this section we show
how this can be done in an explicit way.

Consider the space $\Psi_{\nu,\eta}^\R$
of all real analytic functions
$\psi$ from $(0,\infty)$ to $V_\eta$ satisfying
\begin{eqnarray}
\label{Lewiseqreal}
\eta(T)\psi(x)&=&\psi(x+1)+(x+1)^{-2\nu-1}\eta(ST^{-1})\psi\left(\frac x{x+1}\right)\\
\label{Lewisasym1real}
\psi(x)&=&o(1/x), \quad \text{ as } x\to 0, x>0,\\
\label{Lewisasym2real}
\psi(x)&=&o(1), \quad \text{ as } x\to +\infty, x\in\R.
\end{eqnarray}

Our second main result (see Theorem \ref{periodrealchar})
is a generalization of \cite{Lewis-Zagier2}, Thm.~2, and says
that for  $\Re \nu> -\frac 12$ we have 
$\Psi_{\nu,\eta}^\R
=\{\psi\vert_{(0,\infty)} : \psi\in \Psi_{\nu,\eta}^o\}$.

We summarize the various spaces and mappings considered so far in one diagram:

$$
\begin{diagram} 
\node{\CM_{\nu,\eta}} \arrow{s}
	\node{\CS_{\nu,\eta}} \arrow{w} \arrow{e,t}{L} \arrow{s}
		\node{\Psi_{\nu,\eta}^o} \arrow{s} \arrow{e,t}{\text{res}}
			\node{\Psi_{\nu,\eta}^\R}\\
\node{\tilde \CM_{\nu,\eta}}
	\node{\CF_{\nu,\eta}} \arrow{e}
		\node{\Psi_{\nu,\eta}}\\
	\node[2]{A_{\nu,\eta}^{-\omega}} \arrow{nw,b}{P} \arrow{n,r}{\cong} 		
								\arrow{ne,r}{B}				
\end{diagram}
$$

\bigskip

\section{Automorphic hyperfunctions}\label{Bruggeman}

Let $A$ denote the subgroup of $G$
consisting of diagonal matrices and let $N$ be the
subgroup of upper triangular matrices with $\pm 1$ on the diagonal.
Let $P=AN$ be the group of upper triangular matrices. Finally, let
$K={\rm PSO}(2)={\rm SO}(2)/\{\pm 1\}$ be
the canonical maximal compact subgroup
of $G$. The group $G$ then as a manifold is a direct product
$G=ANK=PK$. For $\nu\in\C$ and
$a=\pm\diag(t,t^{-1})$, $t>0$, let $a^\nu = t^{2\nu}$.
We insert the factor $2$ for compatibility reasons.
Let $(\pi_\nu,V_{\pi_\nu})$ denote the
principal series representation of $G$ with parameter $\nu$. The
representation space $V_{\pi_\nu}$ is the Hilbert space of all functions
$\ph\colon G\to \C$ with
$\ph(anx)=a^{\nu+\frac 12}\ph(x)$ for $a\in A, n\in N, x\in G$,
and $\int_K |\ph(k)|^2\, dk <\infty$
modulo nullfunctions. The representation is
$\pi_\nu(x)\ph(y)=\ph(yx)$.
There is a special vector $\ph_0$ in $V_{\pi_\nu}$ given by
$$
\ph_0(ank)\= a^{\nu+\frac 12}.
$$
This vector is called the \emph{basic spherical function} with parameter
$\nu$. The group $G$ acts on the complex projective line
$\P_1(\C)=\C\cup\{\infty\}$ by linear fractions. This action has
three orbits: the upper half plane $\H^+$, the lower half plane $\H^-$ and the real
projective line $\P_1(\R)=\R\cup\{\infty\}$. The upper half plane can be
identified with $G/K$ via $gK\mapsto g.i$ and $\P_1(\R)$ can be identified
with $P\bs G$ via
$$
P\matrix abcd\ \mapsto\ [c:d].
$$
Our embedding of $\R$ into $\P_1(\R)$ is via $x\mapsto [1:x]$, which can be viewed as
the map
\begin{eqnarray*}
N&\to& P\backslash PwN\\
\matrix 1x01&\mapsto& P\matrix 01{-1}x
\end{eqnarray*}
with the Weyl group element $w=\pm \matrix 01{-1}0$.
Note that $V_{\pi_\nu}$ can also be viewed as a space of sections of a line bundle
over $P\backslash G$. For this bundle the above embedding provides a trivialization
over $\R$. Using the corresponding Bruhat decomposition
$$\pm \matrix abcd= \pm \matrix {c^{-1}}a0c \matrix  01{-1}{\frac dc}$$
for $c\not=0$ we obtain a realization of $V_{\pi_\nu}$ on
$L^2(\R, \frac 1\pi (1+x^2)^{2\nu} dx)$
with the action
$$\pi_\nu\matrix abcdf(x)=(cx-a)^{-2\nu-1} f\left(\frac{dx-b}{cx-a}\right).$$
Transferring the action to $L^2(\R, \frac 1\pi\frac {dx}{1+x^2})$
then yields the action
$$
\pi_\nu\matrix abcd \ph(x) \=
\left(\frac{1+x^2}{(cx-a)^2+(dx-b)^2}\right)^{\nu+\frac
12}\ph\left(\frac{dx-b}{-cx+a}\right)
$$
used in \cite{Brug}.
This is the realization of the principal series we shall work with. Note that in this
realization the basic spherical function is simply the constant function $1$.

Let $\pi_\nu^\omega\subset\pi_\nu^{-\omega}$ be the sets of
analytic vectors and hyperfunction vectors, respectively.
For any  open neighbourhood $U$ of $\P_1(\R)$
inside $\P_1(\C)$ the space $\pi_\nu^{-\omega}$ can be identified with the space
$$
\CO(U\setminus\P_1(\R))/\CO(U),
$$
where $\CO$ denotes the sheaf of holomorphic
functions. This space does not
depend on the choice of $U$. For $U\subseteq\C$ this follows
from Lemma 1.1.2 of
\cite{Schlichtkrull} and generally by subtracting the Laurent series at
infinity. The $G$-action is given by the above formula, where $x$ is replaced
by a complex variable $z$. Note that
any hyperfunction $\al$ on $\P_1(\R)$ has a restriction to $\R$ which
can be represented by a holomorphic function on $\C\setminus\R$.

\begin{proposition}\label{gluingcond} {\rm  (Symmetry of gluing conditions)}
For $f\in\CF_\eta$ the following conditions
are equivalent:
\begin{enumerate}
\item[{\rm(1)}]  $z\mapsto f(z)-z^{-2\nu- 1}\,\eta(S) f\left(\frac{-1}z\right)$
extends holomorphically to $\C\setminus (-\infty,0]$.
\item[{\rm(2)}]  $z\mapsto (1+z^{-2})^{\nu+\frac 12} f\left(\frac{-1}z\right)$ and
$z\mapsto (1+z^2)^{\nu+ \frac 12}\,\eta(S) f(z)$ define the same hyperfunction on
$\R\setminus \{0\}$.
\end{enumerate}
\end{proposition}

\prf
``(2)$\Rightarrow$(1)'' 
Suppose that
$$(1+z^{-2})^{\nu+\frac 12} f\left(\frac{-1}z\right)=(1+z^2)^{-\nu- \frac 12}\,\eta(S)
f(z)+q(z)$$
with a function $q$ that is holomorphic  in a neighborhood of  $\R\setminus \{0\}$.
For $\Re z>0$ we can divide the equation by $(1+z^2)^{\nu+ \frac 12}$ and obtain
$$z^{-2\nu- 1} f\left(\frac{-1}z\right)=\eta(S) f(z)+(1+z^2)^{-\nu- \frac 12}\,q(z).$$
Since $\eta(S)=\eta(S)^{-1}$, this implies the claim.

``(1)$\Rightarrow$(2)''
If (1) holds,  by the same calculation as above we see that for $\Re z>0$ the function
 $$z\mapsto (1+z^2)^{\nu+ \frac 12}\,\eta(S) f(z)-(1+z^{-2})^{\nu+\frac 12}
f\left(\frac{-1}z\right)$$
extends holomorphically to the entire right halfplane. But then the symmetry of this expression
under the
transformation $z\mapsto -\frac 1z$ yields the holomorphic extendability also on the left
halfplane which
proves (2).
\qed

Recall the space 
$A_{\nu,\eta}^{-\omega}=(\pi_\nu^{-\omega}\otimes\eta)^\Ga=
H^0(\Ga,\pi_\nu^{-\omega}\otimes\eta)$
of $\eta$-automorphic hyperfunctions
from (\ref{etaautomhypfcts}).

\begin{proposition}\label{1.1} {\rm (cf. \cite{Brug}, Prop.~2.1)}
There is a bijective linear map
\begin{eqnarray*}
A_{\nu,\eta}^{-\omega} &\to&
\CF_{\nu,\eta}\\
\al &\mapsto & f_\al
\end{eqnarray*}
such that the function $z\mapsto
(1+z^2)^{\nu+\frac12} f_\al(z)$
represents the restriction $\al|_\R$.
\end{proposition}

\prf
The space
$A_{\nu,\eta}^{-\omega}=\left(\pi_\nu^{-\omega}\otimes\eta\right)^\Ga$ can be
viewed as the space of all
$V_\eta$-valued hyperfunctions $\al$ in $\P_1(\R)$
satisfying the invariance condition
$$\pi_\nu(\ga^{-1})\al=\eta(\ga)\al$$
for every  $\ga\in\Ga$. Pick a
representative $f$ for $\al$. The $V_\eta$-valued function
$F\colon z\mapsto (1+z^2)^{-\nu-\frac 12}f(z)$ is holomorphic on
$0<|\Im(z)|<\eps$ for some $\eps>0$.
Note that the invariance of $\al$ under $T$ implies that
for some function $q$, holomorphic on a
neighbourhood of $\R$, we have
\begin{eqnarray*}
\eta(T) f(z)+q(z)
&=&\Big(\pi_\nu(T^{-1}) f\Big)(z)\\
&=&\left(\frac{1+z^2}{1+(z+1)^2}\right)^{\nu+\frac
12} f(z+1)\\
&=&(1+z^2)^{\nu+\frac 12} F(z+1),
\end{eqnarray*}
so that
$$F(z+1)= \eta(T) F(z)+(1+z^2)^{-(\nu+\frac 12)} q(z).$$
Therefore $F$ represents a hyperfunction
on
$\P_1(\R)=\R\cup \{\infty\}$
which is invariant under the translation $z\mapsto z+N$. This
hyperfunction has a representative which is holomorphic in
$\P_1(\C)\setminus\P_1(\R)$.
The freedom in this representative is an
additive constant. So there is a unique
representative $f_\al$ of the form
$$
f_\al(z)\=
\begin{cases}
\ds \frac 12 v_0
+\sum_{k=1}^\infty e^{2\pi i \frac kN z} v_k^+, & z\in\H^+,\\
\ds -\frac 12 v_0
-\sum_{k=1}^\infty e^{-2\pi i \frac kN z} v_k^-, & z\in\H^-.
\end{cases}
$$
So $f_\al\in \CF_\eta$ and $(1+z^2)^{\nu+\frac{1}{2}}f_\al(z)$
represents $\al\vert_\R$.
To show  the injectivity of the map in the Proposition assume that $f_\al=0$. Then $\al$ is supported
in $\{\infty\}$. Since latter set is
not $\Ga$-invariant, $\al$ must be zero.
To see that $f_\al$ lies in $\CF_{\nu,\eta}$, recall
that the invariance of
$\al$ under $S$ implies that
$$
(1+z^{-2})^{\nu+\frac
12}f_\al\left(\frac{1}{-z}\right)
\= (1+z^{2})^{\nu+\frac 12}\eta(S) f_\al\left({z}\right)
+\tilde q(z)
$$
with $\tilde q(z)$ holomorphic on a neighbourhood of $\R\setminus\{0\}$.
Thus Proposition \ref{gluingcond} shows that $f_\al$ satisfies (\ref{Fnueta})
and hence $f_\al\in \CF_{\nu,\eta}$. To finally show surjectivity, let $f\in\CF_{\nu,\eta}$.
Then the function
$$z\mapsto (1+z^2)^{\nu+\frac 12}f(z)$$
represents a hyperfunction $\beta_0$ on $\R$ that satisfies
$\pi_\nu(T^{-1})\beta_0=\eta(T)\beta_0$.
Let $\beta_\infty:=(\pi_\nu\otimes\eta)(S)\beta_0$. Then $\beta_\infty$
is a hyperfunction on $\P_1(\R)\setminus\{ 0\}$ with representative
$z\mapsto (1+z^{-2})^{\nu+\frac 12} \eta(S)f(\frac{-1}z)$.
According to Proposition \ref{gluingcond} the restrictions of  $\beta_0$ and $\beta_\infty$
to  $\P_1(\R)\setminus\{0,\infty\}$ agree.
Thus $\beta_0$ and $\beta_\infty$ are restrictions
of   a  hyperfunction $\beta$ on
$\P_1(\R)$ which is then easily seen to be $S$-invariant.
Using $\beta_0$ we see that the support of
$$(\pi_\nu\otimes\eta)(T)\beta-\beta$$
is contained in $\{\infty\}$.
Using $\beta_\infty$ we see that for $|z|>2$, $z\notin\R$, this
hyperfunction is represented by
\begin{eqnarray*}
&&\hskip -4em
\left(\frac{1+z^2}{1+(z-1)^2}\right)^{\nu+\frac
12}\(1+(z-1)^{-2}\)^{\nu+\frac
12}\eta(TS)
f\(\frac{-1}{z-1}\)\\
&&\qquad-(1+z^{-2})^{\nu+\frac12}
\eta(S)f\(\frac{-1}z\)=\\
&=& (1+z^2)^{\nu+\frac
12}(z-1)^{-2\nu-1}\eta(TS)f\(\frac{-1}{z-1}\)\\
&&\qquad -(1+z^{-2})^{\nu+\frac
12}\eta(S)f\(\frac{-1}{z}\)\\
&=&
(1+z^{-2})^{\nu+\frac 12}\times \\
&&\(
\(\frac
z{z-1}\)^{2\nu+1} \eta(TST^{-1})f\(\frac{z-2}{z-1}\)-\eta\left(
ST^{-1}\right)
f\(\frac{z-1}z\)\)
\end{eqnarray*}
Since $f(z)$ is holomorphic around $z=1$ it follows that
this function is
holomorphic around $z=\infty$. Hence $\beta$ is invariant under $T$.
Now the claim follows because the elements  $S$ and $T$ generate $\Ga$.
\qed

\begin{proposition}\label{1.2}
{\rm (Bruggeman transform; cf. \cite{Brug}, Prop.~2.3)}
For $\al\in A_{\nu,\eta}^{-\omega}$ put
$$
\psi_\al(z)\ := \
f_\al(z)-z^{-2\nu-1}\eta(S)f_\al\left(\frac{-1}z\right),
$$
with $f_\al$ as in Proposition \ref{1.1}. Then the \emph{Bruggeman transform}
$B\colon\al\mapsto\psi_\al$ maps
$A_{\nu,\eta}^{-\omega}$ to $\Psi_{\nu,\eta}$.
It is a bijection if $\nu\notin \frac 12+\Z$.
\end{proposition}

\prf
Let $\al\in A_{\nu,\eta}^{-\omega}$ and define $\psi_\al$ as in the Proposition.
By Proposition
\ref{1.1} the map $\psi_\al$ extends to $\C\setminus (-\infty,0]$.
We compute
\begin{eqnarray*}
&&\hskip -2em
  \psi_\al(z+1)+(z+1)^{-2\nu- 1}\eta(ST^{-1})
     \psi_\al\left(\frac{z}{z+1}\right)= \\
&=&f_\al(z+1)-(z+1)^{-2\nu-1}\eta(S)
   f_\al\left(\frac{-1}{z+1}\right)+ (z+1)^{-2\nu- 1}\eta(ST^{-1})\times \\
& &\quad\ \times \left(f_\al    \left(\frac{z}{z+1}\right)
    -\left(\frac z{z+1} \right)^{-2\nu-1}\eta(S)
    f_\al   \left(\frac{-1}{\frac z{z+1}}\right)     \right).\\
\end{eqnarray*}
Since $\frac z{z+1}=1-\frac 1{z+1}$ and $f_\al\left(1-\frac 1{z+1}\right)=\eta(T)
f_\al\left(\frac{-1}{z+1}\right)$ we see that the two middle summands cancel out.
It remains
\begin{eqnarray*}
&&\hskip -2em\eta(T) f_\al(z)-z^{-2\nu-1}\eta(ST^{-1}
S)f_\al\left(\frac{-z-1}z\right)=\\
&=& \eta(T) f_\al(z)-z^{-2\nu-1}\eta(ST^{-1}ST^{-1})f_\al\left(\frac{-1}z\right)\\
&=& \eta(T)\left(f_\al(z)-z^{-2\nu- 1}\eta(S)f_\al\left(\frac{-1}z\right)\right)\\
&=&\eta(T)\psi_\al(z).
\end{eqnarray*}
Here we have used $ST^{-1}ST^{-1}=TS$. This proves that $\psi_\al$
satisfies the functional equation (\ref{Psinuetaa}).

Next, if
$\nu\in\frac 12+\Z$, then one sees that
$\psi_\al(z)+z^{-2\nu-1}\eta(S)\psi_\al\left(\frac{-1}z\right)$ equals zero and so $\psi_\al$ lies
in
$\Psi_{\nu,\eta}$. If $\nu\notin\frac 12+\Z$ then recall that we take the
standard branch of the logarithm to define $z^{-2\nu-1}$. For $\psi(-1/z)$ one
then takes a complimentary branch and one gets the inversion formula
\begin{equation}\label{brugtrafoinversion}
f_\al(z)\=\frac 1{1+e^{\pm 2\pi  i\nu}}\left(
\psi_\al(z)+z^{-2\nu-1}\eta(S)\psi_\al\left(\frac{-1}z\right)\right)
\end{equation}
for $z\in\H^\pm$. This proves $B\alpha \in \Psi_{\nu,\eta}$ and it only remains to
show that the Bruggeman transform is surjective. But a simple calculation, similar to
the one given above shows that for a holomorphic function
$\psi\colon \C\setminus(-\infty,0]\to V_\eta$ satisfying
(\ref{Psinuetaa}) the function $f\colon \C\setminus \R\to V_\eta$, defined from $\psi$
via the inversion formula (\ref{brugtrafoinversion}), satisfies (\ref{Fetaa}).
If $\psi$ satisfies (\ref{Psinuetab}), then $f$ satisfies (\ref{Fetab})
and (\ref{Fetac}). In view of Proposition \ref{1.1} this, finally, proves the claim.

\qed

\section{Maa\ss\ wave forms}\label{Maass}

Recall the space $\CS_{\nu,\eta}$ of \emph{Maa\ss\ cusp forms} from
(\ref{Snueta}) and consider a $u$ in $\CS_{\nu,\eta}$. Because of
$u(z+N)=u(z)$ the function $u$ has a Fourier series
$$
u(z)\= u(x+iy)\= \sum_{\stackrel{k\in\Z}{k\ne 0}}A_k(y)\, e^{2\pi i\frac kN
x}\, v_k
$$
for some $v_k\in V_\eta$. The differential equation $\Delta u=(\frac 14
-\nu^2)u$ implies a differential equation for $A_k(y)$ which implies that it
must be a linear combination of $I$ and $K$-Bessel functions.
The fact that $u$ is square integrable rules out the $I$-Bessel functions, so
$$
A_k(y)\= \sqrt y K_\nu\left( 2\pi \frac{|k|}N y\right)
$$
times a constant which we can assume to be $1$ by multiplying it to $v_k$. 
By Theorem 3.2 of \cite{Iwaniec} it follows that the norms $\norm{v_k}$ are bounded as $|k|\to\infty$.
The functional equation
$u(z+1)=\eta(T)u(z)$ is reflected in the fact that the $v_k$ are eigenvectors
of $\eta(T)$, since we get $\eta(T)v_k=e^{2\pi i\frac kN}v_k$. Now set
\begin{equation}\label{fudef}
\ f_u(z)\ :=\  \begin{cases}\ds\sum_{k>0}k^\nu e^{2\pi i\frac kN z}v_k, & \Im(z)>0,\\
\ds -\sum_{k<0}|k|^\nu e^{2\pi i\frac kN z}v_k,& \Im(z)<0.
\end{cases}
\end{equation}
From the construction it is clear that $f_u$ satisfies (\ref{Fetaa}) - (\ref{Fetac}),
i.e. $f_u\in \CF_\eta$. It will play the role of our earlier $f_\al$ (cf. Proposition \ref{1.1}), so
we define
\begin{equation}\label{psiudef}
\psi_u(z)\ :=\ f_u(z)-z^{-2\nu- 1}\eta(S) f_u\(\frac{-1}z\).
\end{equation}

\begin{lemma}\label{Lewistrafolemma}
For  $\Re \nu>-\frac 12$ the equation (\ref{fudef})  and (\ref{psiudef})  define   linear maps
$$
\begin{array}{rcl}
\CS_{\nu,\eta}&\to& \CF_{\nu,\eta}\\
 u&\mapsto& f_u
\end{array}
\quad\text{and}\quad
\begin{array}{rcl}
\CS_{\nu,\eta}&\to& \Psi_{\nu,\eta}^o.\\
 u&\mapsto& \psi_u
\end{array}
$$
\end{lemma}

\prf
To prove that $f_u\in \CF_{\nu,\eta}$  we will need the
following two Dirichlet series. For $\eps =0,1$ set
$$
L_\eps(u,s)\ :=\  \sum_{k\ne 0} {\rm sign}(k)^\eps \left(\frac N{|k|}\right)^s
v_k.
$$
We will relate $L_0$ and $L_1$ to $u$ by the Mellin transform. For this let
\begin{equation}\label{uepsilon}
u_0(y)\=\frac 1{\sqrt y} u(iy),\qquad u_1(y)\=\frac{\sqrt y}{2\pi i}u_x(iy),
\end{equation}
where $u_x=\frac{\partial}{\partial x}u$.
Next define
\begin{equation}\label{Ldachepsilon}
\hat L_\eps(u,s)\ :=\ \int_0^\infty u_\eps(y) y^s\frac{dy}y.
\end{equation}
Plugging in the Fourier series of $u$ and using the fact that
$$
\int_0^\infty K_\nu(2\pi y)y^s\frac{dy}y
\=\Ga_\nu(s)
\ :=\ \frac 1{4\pi^s}\Ga\left(\frac{s-\nu}2\right)\Ga\left(\frac{s+\nu}2\right),
$$
we get
$$
\hat L_0(u,s)\=\Ga_\nu(s) L_0(u,s),
$$
and similarly,
$$
\hat L_1(u,s)\=\Ga_\nu(s+1) L_1(u,s).
$$
On the other hand, the usual process of splitting the Mellin integral and
using the functional equations
$$
u_\eps\left(\frac 1y\right)\= (-1)^\eps y\,\eta(S)\, u_\eps(y),\qquad
\eps=0,1,
$$
(which can be checked using the Taylor series of $u$), one gets that $\hat L_\eps$
extends to an entire function and satisfies the
functional equation,
$$
\hat L_\eps(u,s)\= (-1)^\eps\,\eta(S)\, \hat L_\eps(u,1-s).
$$
With a similar, even easier computation one gets
$$
\int_0^\infty y^s \left(f_u(iy)-(-1)^\eps f_u(-iy)\right)\, \frac {dy}y
\=\frac{\Ga(s)N^\nu}{(2\pi)^s}L_\eps(u,s-\nu).
$$
This implies that the {\emph Mellin transforms}
$M^\pm f(s):= \int_0^\infty y^s f(\pm iy)\,\frac{dy}y$ can be calculated as
\begin{eqnarray*}
&&\hskip -2em M^\pm f(s)=
   {\textstyle \pm\frac{\Ga(s) N^\nu}{2(2\pi)^s}\left(L_0(u,s-\nu)\pm
    L_1(u,s-\nu)\right)}\\
&=& {\textstyle \pm N^\nu \pi^{-\nu-\frac 32}
    \Ga\left(\frac{s+1}2\right)\Ga\(\frac{2\nu+2-s}{2}\)\sin\pi\(\nu+1-\frac s2\)\,
    \hat L_0(u,s-\nu)}\\
&&{\textstyle +N^\nu\pi^{-\nu-\frac 12}\Ga\(\frac s2\)\Ga\(\frac{2\nu+1-s}2\)
     \sin\pi\(\nu+\frac 12-\frac s2\)\,\hat L_1(u,s-\nu)}.
\end{eqnarray*}
The last identity follows from the standard equations
$$
{\textstyle
\Ga\(\frac x2\)\Ga\(\frac{x+1}2\)\= \Ga(x) 2^{1-x}\sqrt\pi,\qquad
\Ga(x)\Ga(1-x)\=\frac\pi{\sin\pi x}.}
$$
Thus the Mellin transform $M^\pm f(s)$ is seen to be holomorphic for
$\Re(s)>0$ and rapidly decreasing on any vertical strip. The Mellin
inversion formula yields for $C>0$,
$$
f_u(\pm iy)\=\frac 1{2\pi i}\int_{\Re(s)=C} y^{-s} M^\pm f_u(s)\, ds.
$$
This extends to any $z\in\C\setminus \R$ to give
$$
f_u(z)\=\frac 1{2\pi i} \int_{\Re(s)=C} e^{\pm\frac{\pi}2is}z^{-s}\,
M^\pm f_u(s)\, ds
$$
for $z\in\H^\pm$. For $0<C<2\Re \nu+1$ (here we need $\Re \nu>-\frac 12$)
it follows that
$$
\psi_u(z)\ =
 \frac 1{2\pi i}
\int_{\Re(s)=C}\left(e^{\pm \frac\pi 2 is}z^{-s} - e^{\mp\frac \pi
2is}z^{-2\nu-1}z^s\eta(S)\right) M^\pm f_u(s)\, ds.
$$
Writing this as the difference of two integrals, substituting $s$ in the
second integral with $2\nu +1-s$ and shifting the contour we arrive at
the formula
\begin{equation}\label{intformforpsiu}
\begin{array}{rl}{\frac 1{2\pi i} \int_{\Re(s)=C}}&
\Big(e^{\pm\frac \pi 2 is}z^{-s} M^\pm f_u(s)\\
&-e^{\mp\frac\pi 2i(2\nu+1-s)} z^{-s}\eta(S)M^\pm f_u(2\nu+1-s)\Big)\, ds.
\end{array}
\end{equation}
for $\psi_u$. Using the identities
\begin{eqnarray*}
\pm e^{\pm\frac\pi 2is}\cos\pi\(\nu +\frac 12 -\frac s2\)\mp e^{\mp\frac\pi
2i(2\nu+1-s)}\cos \pi\frac s2 &=& i\sin\pi\(\nu+\frac 12\),\\
 e^{\pm\frac\pi 2is}\sin\pi\(\nu +\frac 12 -\frac s2\)+ e^{\mp\frac\pi
2i(2\nu+1-s)}\sin \pi\frac s2 &=& \sin\pi\(\nu +\frac 12\),
\end{eqnarray*}
and the functional equation of $\hat L_\eps$ we see that the integrand
of (\ref{intformforpsiu}) equals
\begin{equation}\label{intformforpsiu2}
\begin{array}{rl}
z^{-s}N^\nu \sin\pi\(\nu+\frac 12\) &
\left[  \pi^{-\nu-\frac 32} \Ga\(\frac{s+1}2\)\Ga\(\frac{2\nu+2-s}2\) i\, \hat
L_0(u,s-\nu)+\right.\\
&\quad+\left.  \pi^{-\nu-\frac 12} \Ga\(\frac s2\)\Ga\(\frac{2\nu+1-s}2\)\,
\hat L_1(u,s-\nu)\right].
\end{array}
\end{equation}
Since this expression is independent of whether $z$ lies in $\H^+$ or $\H^-$,
it follows that $f_u(z)-z^{-2\nu-1}\eta(S) f_u\(\frac{-1}z\)$ extends to a
holomorphic function on $\C\setminus (-\infty,0]$, i.e., the function $f_u$
indeed lies in the space $\CF_{\nu,\eta}$. The linearity of the map is clear.

It remains to show that $\psi_u\in \Psi_{\nu,\eta}^o$ .
Note that in view of   $f_u\in \CF_{\nu,\eta}$ Proposition \ref{1.1} shows that
the function $z\mapsto (1+z^2)^{\nu+\frac 12}f_u(z)$
represents a hyperfunction $\al_u\in A_{\nu,\eta}^{-\omega}$. Then,
according to Proposition \ref{1.2} we have $\psi_u=B(\alpha_u)$ so
$\psi_u$ satisfies (\ref{Psinuetab}).
The asymptotic property  (\ref{Lewisasymcomplex})
now follows from the integral representation (\ref{intformforpsiu})
with the $C$ chosen there. More precisely, the bound $O(|z|^{-C})$ follows directly from
(\ref{intformforpsiu2}) since the integrant divided by $z^{-s} $ is of $\pi$-exponential decay,
whereas
the $O(1)$-bound is obtained by moving the contour slightly to the left of the imaginary
axis picking up the residue at $0$ which is proportional to $1$
(see \cite{Lewis-Zagier2}, \S I.4 for more details on this type of argument).
\qed

\begin{lemma}\label{Poissonexp}
For $0\not=k\in \Z$ let $\alpha_k$ be the hyperfunction on $\P_1(\R)$ represented by
$(1+z^2)^{\nu+\frac 12} f_k(z)$ with
$$f_k(z)=\begin{cases}
\sign(k)\cdot e^{2\pi i\frac{k}{N}z}&\text{ for } \ \sign(k)\cdot\Im(z)>0\\
0& \text{ for } \ \sign(k)\cdot\Im(z)<0.\\
\end{cases}
$$
Then we have that
$$ \sp{\pi_{-\nu}\matrix{\sqrt b}{\frac a{\sqrt b}}0{\frac 1{\sqrt b}}\ph_0,\al_k}
$$
equals 
$$
2\,\sign(k)\left(\frac N{|k|}\right)^{\nu}
\frac{\pi^{-\nu-\frac 12}}{\Gamma(\frac 12-\nu)}\sqrt{b}\
K_\nu\left(2\pi \frac{|k|}{N} b\right)
e^{2\pi i \frac{k}{N} a},
$$
where is $K_\nu$ the $K$-Bessel function with parameter $\nu$.
\end{lemma}

\prf
For this we will need the following identity (cf. \cite{Brug}, \S4 or
\cite{Ter}, p.136)
\begin{equation}\label{Brugidentitaet}
\int_{-\infty}^\infty
\left(\frac 1{y^2+(\tau-x)^2}\right)^{\frac 12 -\nu}
  e^{2\pi i k\tau}\ dt
=
\frac{2\pi^{\frac 12 -\nu}|k|^{-\nu} }{\Gamma\left(\frac 12-\nu\right)}
y^\nu K_\nu(2\pi|k|y) e^{2\pi ikx}.
\end{equation}
Note that $g=\matrix{\sqrt b}{\frac a{\sqrt b}}0{\frac 1{\sqrt b}}$  satisfies
$g\cdot i=a+ib$. Therefore, by abuse of notation, we write
$P(\alpha_k)(a+ib)$ for
$\sp{\pi_{-\nu}\matrix{\sqrt b}{\frac a{\sqrt b}}0{\frac 1{\sqrt b}}\ph_0,\al_k}$.
According to \cite{Brug}, \S4, we can calculate
\begin{eqnarray*}
P(\alpha_k)(a+ib)
&=&\Big\langle \left(\frac{1+x^2}{b+\left(\frac x{\sqrt b} - \frac a{\sqrt b}\right)^2}
                \right)^{-\nu+\frac 12}, \alpha_k\Big\rangle\\
&=&b^{-\nu+\frac 12}
   \Big\langle \left(\frac{1+x^2}{b^2+\left(x - a\right)^2}
                \right)^{-\nu+\frac 12}, (1+z^2)^{\nu+\frac 12} f_k(z)\Big\rangle\\
&=&\sign(k)b^{-\nu+\frac 12}\frac 1\pi \int_{-\infty}^\infty
    \left(\frac{1}{b^2+\left(x - a\right)^2}\right)^{-\nu+\frac 12}e^{2\pi i\frac kN x}\, dx\\
&=&2\,\sign(k)\left(\frac N{|k|}\right)^{\nu}
\frac{\pi^{-\nu-\frac 12}}{\Gamma(\frac 12-\nu)}\sqrt{b}\
K_\nu\left(2\pi\frac{|k|}{N} b\right)
e^{2\pi i \frac{k}{N} a},
\end{eqnarray*}
where in the last step we have used (\ref{Brugidentitaet}).
\qed

\begin{theorem}\label{Lewistrafobijektiv}
{\rm (Lewis transform; cf. \cite{Lewis-Zagier2}, Thm.~1.1)}
For $\nu\not\in \frac 12+\Z$ and $\Re \nu> -\frac 12$
the Lewis transform is a bijective linear map from
the space of Maa\ss\ cusp
forms $\CS_{\nu,\eta}$ to the space $\Psi_{\nu,\eta}^o$ of period functions.
\end{theorem}

\prf
We begin by showing that the Lewis transform is injective on $\CS_{\nu,\eta}$.
This  will be done by proving that we can recover
$u$ from $\psi_u$, where we use the notation from Lemma \ref{Lewistrafolemma}.
The hypothesis $\nu\not\in \frac 12+\Z$ guarantees that we can recover
$f_u$ from $\psi_u$ via a simple algebraic manipulation (cf. the proof of
Proposition \ref{1.2}).
Thus it suffices to express $u$ in terms of $\alpha_u$. But applying
Lemma \ref{Poissonexp} to the summands in the defining formula (\ref{fudef})
for $f_u$, we obtain
\begin{equation}\label{Poissontransformula}
\begin{array}{rcl}
P\alpha_u(a+ib)
&=&\frac{2\, N^\nu}{\pi^{\nu+\frac 12}\Gamma(\frac 12 -\nu)}\sqrt b\, \sum_{k\not=0}
    K_\nu\left(2\pi\frac{|k|}{N} b\right) e^{2\pi i \frac kN a} v_k\\
&&\\
&=&\frac{2\, N^\nu}{\pi^{\nu+\frac 12}\Gamma(\frac 12 -\nu)}u(a+ib).
\end{array}
\end{equation}

\bigskip

It remains to show that  $L(\CS_{\nu,\eta})=\Psi_{\nu,\eta}^o$.
To do this pick $\psi\in \Psi_{\nu,\eta}^o$.
According to Propositions \ref{1.1} and  \ref{1.2}
we can find
a hyperfunction $\alpha\in A_{\nu,\eta}^\Gamma$ represented by the
function $(1+z^2)^{\nu+\frac 12} f$ with $f\in \CF_{\nu,\eta}$ such that
\begin{eqnarray*}
\psi(z)&=&f(z)-z^{-2\nu-1} \eta(S) f\left(-\frac 1z\right),\\
f(z)&=&\frac 1{1+e^{\pm 2\pi i\nu}}\left( \psi(z)+z^{-2\nu- 1}\eta(S)
        \psi\left(\frac{-1}z\right)\right)
\end{eqnarray*}
for $z\in \H^\pm$. The function $f$ admits a Fourier expansion of the form
$$
f(z)\=
\begin{cases}
\ds \frac 12 v_0
+\sum_{k=1}^\infty e^{2\pi i \frac kN z} v_k, & z\in\H^+,\\
\ds -\frac 12 v_0
-\sum_{k=1}^\infty e^{-2\pi i \frac kN z} v_{-k}, & z\in\H^-.
\end{cases}
$$
The asymptotic property (\ref{Lewisasymcomplex})
of $\psi$ implies that
\begin{eqnarray*}
\psi(z)&=&O(|z|^{-C})\\
z^{-2\nu-1}\eta(S){\textstyle \psi\left(-\frac 1z\right)} &=&O(|z|^{-2\nu-1})
\end{eqnarray*}
for $\z\in \C\setminus (-\infty,0]$. Since $2\nu+1>0$ this implies that there is a
constant $\epsilon>0$ such that
$$f(x+i y)=O(|y|^{-\eps})$$
locally uniformly in $x$.
Since $f$ is periodic, this shows $v_0=0$.
Note that $K_\nu\left(t\right) \sim   e^{-t}\sqrt{\frac \pi {2t}}$.
Therefore we have
\begin{equation}\label{Akyasymptotik}
A_k(y)\= \sqrt y K_\nu\left( 2\pi \frac{|k|}N y\right)\sim
e^{-2\pi \frac{|k|}N y} \sqrt{\frac N {4\, |k|}}
\end{equation}
uniformly in $k$, which implies that
$$
u(z)\ := \ u(x+iy)\ :=\ \sum_{\stackrel{k\in\Z}{k\ne 0}}A_k(y)\, e^{2\pi i\frac kN
x}\, v_k
$$
defines a smooth function on  $\H^+$. Taking the derivatives termwise, we see that
$u$ satisfies (\ref{Mnuetac}), i.e. is contained in the range of the Poisson transform.
Now Lemma \ref{Poissonexp} shows that
$$\sp{\pi_{-\nu}\matrix{\sqrt b}{\frac a{\sqrt b}}0{\frac 1{\sqrt b}}\ph_0,\al}=u$$
and (\ref{Poissonisom}) implies  $u\in \tilde \CM_{\nu,\eta}$. Note that
(\ref{Snueta}) is a consequence of $v_0=0$. Thus in order to show that
$\psi\in L(\CS_{\nu,\eta})$, it only remains to show that $u$ satisfies (\ref{Mnuetab}).
But (\ref{Akyasymptotik}) implies that $u$ rapidly decreases towards the cusp and hence
the finite volume of the fundamental domain proves the square integrability of $u$.
\qed

As a consequence of this proof we see that for $\eta$ the trivial representation,
our Lewis transform coincides with $\frac 12 \pi^{\nu+\frac 12}\Ga\(\frac 12 -\nu\)$
times the one given in \cite{Lewis-Zagier2}.

\section{Characterizing period functions on $\R^+$}
Let $T'=\left(\begin{array}{cc}1 & 0 \\1 & 1\end{array}\right)=\( STS^{-1}\)^{-1}$.

\begin{lemma}\label{realasymptoticslemma}
{\rm (cf. \cite{Lewis-Zagier2}, \S III.3)}
If a smooth function $\psi\colon (0,\infty)\to V_\eta$ satisfies
(\ref{Lewiseqreal}) with $\nu\not\in \frac 12+\Z$, then it has the
following asymptotic expansions:
\begin{eqnarray*}
\psi(x)&\underset{x\to0}{\sim}&
x^{-2\nu-1} Q_0{\textstyle\left(\frac 1x\right)}+\sum_{m=-1}^\infty C_m^* x^m,\\
\psi(x)&\underset{x\to\infty}{\sim}&
x^{-2\nu-1} Q_\infty{\textstyle\left(\frac 1x\right)}
  +\sum_{m=-1}^\infty (-1)^mC_m^* x^{-m-2\nu-1},\\
\end{eqnarray*}
where the $Q_0,Q_\infty\colon R\to \C$ are smooth functions with
\begin{eqnarray*}
Q_0(x+1)&=&\eta(T')Q_0(x),\\
Q_\infty(x+1)&=&\eta(T)Q_\infty(x),
\end{eqnarray*}
and the $C_m^*$ can be calculated from the Taylor coefficients
$C_m:=\frac 1{m!}\psi^{(m)}(1)\in V_\eta$ of $\psi$ in $1$ via
\begin{equation}\label{Cstar}
C_m^*=\frac 1 {m+2\nu+1} \sum_{k=0}^M(-1)^m B_k\,
\begin{pmatrix} m+2\nu+1\\ k\end{pmatrix} C_{m-1-k}.
\end{equation}
Here the $B_k$ are the Bernoulli numbers. If $\psi$ is real analytic, then so are
$Q_0$ and $Q_\infty$.
\end{lemma}

\prf
For $\Re \nu>0$ set
\begin{eqnarray*}
Q_0(x)
&:=& x^{-2\nu-1} \psi\left(\frac 1x\right)
-\sum_{n=0}^\infty (n+x)^{-2\nu-1} \eta\big(T(T')^n\big)^{-1}
\psi\left(1+\frac 1{n+x}\right)\\
\text{and}&&\\
Q_\infty(x)
&:=& \psi\left(x\right)- \sum_{n=1}^\infty (n+x)^{-2\nu-1} \eta\big(T'T^{n-1}\big)^{-1}
    \psi\left(1-\frac 1{n+x}\right).
\end{eqnarray*}
Then we have
\begin{eqnarray*}
&&\hskip -2em Q_0(x+1)-\eta(T') Q_0(x)=\\
&=&(x+1)^{-2\nu-1}\psi\left(\frac 1{x+1}\right)-x^{-2\nu-1}\eta(T')\psi\left(\frac
1{x}\right)\\
& &-\sum_{n=0}^\infty (n+1+x)^{-2\nu-1}\eta(T(T')^n)^{-1} \psi\left(1+\frac
1{n+1+x}\right)\\
& &+\sum_{n=0}^\infty (n+x)^{-2\nu-1}\eta(T')\eta(T(T')^n)^{-1}
    \psi\left(1+\frac 1{n+x}\right)\\
&=&(x+1)^{-2\nu-1}\psi\left(\frac 1{x+1}\right)-x^{-2\nu-1}\eta(T')\psi\left(\frac
1{x}\right)\\
& &-\sum_{n=1}^\infty (n+x)^{-2\nu-1}\eta(T(T')^{n-1})^{-1} \psi\left(1+\frac
1{n+x}\right)\\
& &+\sum_{n=0}^\infty (n+x)^{-2\nu-1}\eta(T(T')^{n-1})^{-1}
    \psi\left(1+\frac 1{n+x}\right)\\
&=&(x+1)^{-2\nu-1}\psi\left(\frac 1{x+1}\right)-x^{-2\nu-1}\eta(T')\psi\left(\frac
1{x}\right)\\
& &+x^{-2\nu-1}\eta(T(T')^{-1})^{-1}
   \left(\eta(T)\psi\left(\frac 1{x}\right)\right.\\
& &\left.   -\left(1+\frac 1x\right)^{-2\nu-1}\eta(ST^{-1})
       \psi\left(\frac{\frac 1x}{1+\frac 1x}\right)\right)\\
&=&0,
\end{eqnarray*}
since $T^{-1}ST^{-1}=(T')^{-1}$. Similarly we calculate
\begin{eqnarray*}
&&\hskip -3em Q_\infty(x+1)-\eta(T) Q_\infty(x)=\\
&=&\psi\left({x+1}\right)-\eta(T)\psi\left({x}\right)\\
& &-\sum_{n=1}^\infty (n+1+x)^{-2\nu-1}\eta(T'T^{n-1})^{-1} \psi\left(1-\frac
1{n+1+x}\right)\\
& &+\sum_{n=1}^\infty (n+x)^{-2\nu-1}\eta(T)\eta(T'T^{n-1})^{-1}
    \psi\left(1-\frac 1{n+x}\right)\\
&=&\psi\left({x+1}\right)-\eta(T)\psi\left({x}\right)\\
& &-\sum_{n=2}^\infty (n+x)^{-2\nu-1}\eta(T'T^{n-2})^{-1} \psi\left(1-\frac 1{n+x}\right)\\
& &+\sum_{n=1}^\infty (n+x)^{-2\nu-1}\eta(T'T^{n-2})^{-1}
    \psi\left(1-\frac 1{n+x}\right)\\
&=&\psi\left({x+1}\right)-\eta(T)\psi\left({x}\right)\\
& &+(1+x)^{-2\nu-1}\eta(T'T^{-1})^{-1}
   \psi\left(1-\frac 1{x+1}\right)\\
&=&0.
\end{eqnarray*}

For general $\nu$ we write
\begin{eqnarray*}
Q_0(x)
&:=& x^{-2\nu-1} \psi\left(\frac 1x\right)- \sum_{m=0}^M C_m\zeta(m+2\nu+1,x)\\
&&-\sum_{n=0}^\infty (n+x)^{-2\nu-1} \eta\big(T(T')^n\big)^{-1}
    \left(\psi\left(1+\frac 1{n+x}\right)-\sum_{m=0}^M \frac{C_m}{(n+x)^m}\right)\\
Q_\infty(x)
&:=& \psi\left(x\right)- \sum_{m=0}^M (-1)^m C_m\zeta(m+2\nu+1,x+1)\\
&&-\sum_{n=1}^\infty (n+x)^{-2\nu-1} \eta\big(T'T^n\big)^{-1}
    \left(\psi\left(1-\frac 1{n+x}\right)-\sum_{m=0}^M \frac{C_m}{(n+x)^m}\right)\\
\end{eqnarray*}
with the Hurwitz zeta function
$\zeta(a,x):=\sum_{n=0}^\infty\frac 1{(n+x)^{a}}$. Since the Hurwitz zeta function
satisfies
\begin{equation}\label{hurwitzzetaasymptotics}
\zeta(a,x)\underset{x\to\infty}{\sim} \frac{1}{a-1}
\sum_{k\ge 0}(-1)^k B_k \ \begin{pmatrix}k+a-2\\ k\end{pmatrix}\ x^{-a-k+1}
\end{equation}
we find
\begin{eqnarray*}
\psi(x)
&=& x^{-2\nu-1}Q_0(x^{-1})+\sum_{m=0}^M C_m\zeta(m+2\nu+1,x^{-1})
x^{-2\nu-1}\\
& & +\underbrace{\sum_{n=0}^\infty (x^{-1}+n)^{-2\nu-1}
    \left(\psi\left(1+\frac1{n+x^{-1}}\right)
    -\sum_{m=0}^M \frac{C_m}{\left(n+x^{-1}\right)^n}\right).}_{=O(x^{2\nu+1+M})}
\end{eqnarray*}
From this one derives the asymptotics  for $x\to 0$ using (\ref{hurwitzzetaasymptotics}),
see \cite{Lewis-Zagier2},  \S III.3 for details.
The asymptotics for $x\to\infty$ is shown analogously and the last claim is obvious.
\qed

\begin{remark}\label{transferop}
{\rm
\begin{enumerate}
\item[(i)]
If $\psi(x)=o(x^{\min(1,2\Re \nu+1)})$ for $x\to0$, then
$Q_0=0$ by periodicity,  i.e., $\psi$ is an eigenfunction for the \emph{transfer operator}
$$\mathcal L_0\psi(x):=
x^{-2\nu-1}\sum_{n=0}^\infty (n+x^{-1})^{-2\nu-1} \eta\big(T(T')^n\big)^{-1}
    {\textstyle \psi\left(1+\frac 1{n+x^{-1}}\right)}.$$
Moreover we have $C_{-1}^*=0$.
\item[(ii)] If  $\psi(x)=o(x^{\min(0,\Re \nu)})$ for $x\to\infty$, then $Q_\infty=0$,
i.e., $\psi$ is an eigenfunction for the \emph{transfer operator}
$$\mathcal L_\infty\psi(x):=
\sum_{n=1}^\infty (n+x)^{-2\nu-1} \eta\big(T'T^n\big)^{-1}
   {\textstyle \psi\left(1+\frac 1{n+x}\right)}.$$
Moreover we have $C_{-1}^*=0$.
\item[(iii)] If $C_{-1}^*=0$, then $C_0=0$,
and if $Q_0=Q_\infty=0$ we have the equations
\begin{eqnarray}\label{psitransfereqQzero}
\psi(x)&=&
x^{-2\nu-1}\sum_{n=0}^\infty (n+x^{-1})^{-2\nu-1} \eta\big(T(T')^n\big)^{-1}
  {\textstyle  \psi\left(1+\frac 1{n+x^{-1}}\right)}\\
\label{psitransfereqQinfty}
\psi(x)&=&\sum_{n=1}^\infty (n+x)^{-2\nu-1} \eta\big(T'T^n\big)^{-1}
         {\textstyle \psi\left(1+\frac 1{n+x}\right)}.
\end{eqnarray}
In this case, we can analytically extend $\psi$ to $\C\setminus(-\infty,0]$ via
$$\psi(z):=\sum_{\gamma\in Q_n} (\psi\vert_{\nu,\eta} \gamma)(z),$$
where $Q$ is the semigroup generated by $T$ and $T'$,
$Q_n$ is the set of $T$--$T'$-words of length $n$ in $Q$, and
\begin{equation}\label{slashwirkung}
(\psi\vert_{\nu,\eta} \gamma)(z):=(cz+d)^{-2\nu-1}
\eta(\gamma)^{-1}\psi(\gamma\cdot z)
\end{equation}
is a well defined right semigroup action (cf. \cite{HMM03}, \S~3,
and \cite{Lewis-Zagier2},  \S III.3).
The analytically continued function $\psi$ still satisfies
(\ref{psitransfereqQzero}) and (\ref{psitransfereqQinfty}).
Therefore we can mimick the proof of Lemma
\ref{realasymptoticslemma} and use the Taylor expansion in $1$ to find
\begin{equation}\label{cplxasymp0viahurwitz}
\psi(z)=
\sum_{m=1}^M
C_m\zeta(m+2\nu+1,z^{-1})z^{-2\nu-1}+O\big(|\zeta(2\nu+M+2,z^{-1})|\big)
\end{equation}
for $|z|\to 0$ and
\begin{equation}\label{cplxasympinftyviahurwitz}
\psi(z)=
\sum_{m=1}^M (-1)^m C_m\zeta(m+2\nu+1,z+1)+O\big(|\zeta(2\nu+M+2,z)|\big),\\
\end{equation}
for  $|z|\to\infty$.
Now we use the following version of (\ref{hurwitzzetaasymptotics})
which can be found in \cite{EMOT}, \S~1.18:
\begin{equation}\label{cplxzetaasymptotics}
\begin{array}{rcl}
\zeta(a,z)&=&z^{1-a}\frac{\Gamma(a-1)}{\Gamma(a)}+ \frac 12 z^{-a}
             +\sum_{n=1}^N B_{2n} \frac{\Gamma(a+2n-1)}{\Gamma(a)(2n)!}z^{1-2n-a}\\
&&+ O(|z|^{-2N-1-a})^{\phantom{\Big(\big)}}\\
\end{array}
\end{equation}
for $\Re a>1$ and $z\in \C\setminus (-\infty,0]$. Then
(\ref{cplxasymp0viahurwitz}) and (\ref{cplxasympinftyviahurwitz})
result in
\begin{equation}\label{cplxasymp0viahurwitz'}
\psi(z)=O(1) \quad \text{for }|z|\to 0
\end{equation}
 and
\begin{equation}\label{cplxasympinftyviahurwitz'}
\psi(z)=O(|z|^{-2\nu-1})\quad \text{for }|z|\to \infty.
\end{equation}
\qed
\end{enumerate}
}\end{remark}

\begin{remark}
{\rm
One can use the slash action (\ref{slashwirkung}) to rewrite the
real version (\ref{Lewiseqreal}) of the Lewis equation in the form
$$\psi=\psi\vert_{\nu,\eta} T+\psi\vert_{\nu,\eta} T'.$$
\qed
}
\end{remark}

\begin{theorem}\label{periodrealchar}
{\rm (cf. \cite{Lewis-Zagier2}, Thm.~2)} Suppose that $\Re \nu> -\frac 12$. Then
$$\Psi_{\nu,\eta}^\R=\{\psi\vert_{(0,\infty)} : \psi\in \Psi_{\nu,\eta}^o\}.$$
\end{theorem}

\prf
Note first that  property (\ref{Lewisasymcomplex}) of $\psi\in \Psi_{\nu,\eta}^o$
trivially implies (\ref{Lewisasym1real}) and  (\ref{Lewisasym2real}) for
$\psi\vert_{(0,\infty)}$. Therefore it only remains to show that
each element of  $\Psi_{\nu,\eta}^\R$ occurs as the restriction of some
$\psi\in \Psi_{\nu,\eta}^o$. To this end we fix a $\tilde \psi\in \Psi_{\nu,\eta}^\R$.
Since (\ref{Lewisasym1real}) and  (\ref{Lewisasym2real}) hold for
$\tilde \psi$, we can apply  Remark \ref{transferop} to it. Thus $\tilde \psi$
has an analytic continuation to $\C\setminus (-\infty,0]$ (still denoted by $\tilde \psi$)
and the asymptotics (\ref{cplxasymp0viahurwitz'}) and (\ref{cplxasympinftyviahurwitz'})
shows that $\tilde \psi$ indeed satisfies (\ref{Lewisasymcomplex}).
\qed

\section{A Converse Theorem}

\begin{theorem}
Let $v_k\in V_\eta$ for $k\in\Z\setminus\{0\}$ such that $\eta(T) v_k=e^{2\pi
i\frac kN}v_k$ and that the two Dirichlet series
$$
L_\eps(s)\=\sum_{k\ne 0} {\rm sgn}(k)^\eps\(\frac N{|k|}\)^s \, v_k,\qquad
\eps=0,1
$$
converge for $\Re(s)>>0$. Assume that $\hat L_\eps(s)=\Ga_\nu(s+\eps)
L_\eps(s)$ extends to an entire function with
$$
\hat L_\eps(s)\= (-1)^\eps\,\eta(S)\, \hat L_\eps(1-s).
$$
Then the function $u$ given by
$$
u(z)\=\sum_{k\ne 0}\sqrt y K_\nu\( 2\pi\frac{|k|}Ny\) \, e^{2\pi i\frac kN
x}\, v_k
$$
lies in $\CS_{\nu,\eta}$.
\end{theorem}

\prf
The Dirichlet series give rise to an inverse Mellin transform $f$ as in
Section \ref{Maass}. Now follow the argumentation of that section.
\qed

{\small
\noindent
University of Exeter,  Mathematics, Exeter EX4
4QE, England.\\
{\tt a.h.j.deitmar@ex.ac.uk}

\noindent
Universit\"at Paderborn, Fakult\"at f\"ur Elektrotechnik, Informatik und
Mathematik, Warburger Str.  100, 33098 Paderborn, Germany.\\
\tt hilgert@math.uni-paderborn.de}

\end{document}